\documentclass[conference]{IEEEtran}
\IEEEoverridecommandlockouts
\usepackage{cite}
\usepackage{amsmath,amssymb,amsfonts}
\usepackage{algorithmic}
\usepackage{graphicx}
\usepackage{textcomp}
\usepackage{xcolor}
\def\BibTeX{{\rm B\kern-.05em{\sc i\kern-.025em b}\kern-.08em
    T\kern-.1667em\lower.7ex\hbox{E}\kern-.125emX}}
\begin{document}

\title{An Overview of Application of Optimization Models Under Uncertainty to the Unit Commitment Problem \\
}

\author{
	\IEEEauthorblockN{
		Zambrano, Angel
	}
	\IEEEauthorblockA{
		\textit{Department of Electrical Engineering and Computer Science} \\
		\textit{Technical University of Berlin}\\
		Berlin, Germany \\
		guerrazambrano@campus.tu-berlin.de
	}
}

\maketitle

\begin{abstract}
Optimization models have been broadly used withinside the energy industry as useful decision-making systems for scheduling and dispatching electric powered energy resources; this is applied in a system called unit commitment (UC).

Unit Commitment seeks the maximum price adequate generator commitment scheme for an electric powered energy device to satisfy a certain demand, at the same time as fulfilling the operational constraints on transmission models and computational resources. Taking into account risk variability in those processes and checking out their comparative overall performance for a single or multi-stage energy model as a function of monetary performance in addition to the risk related to the commitment decisions.

Stochastic programming and Robust optimization are by a vast majority the most widely studied methodologies for UC under net load uncertainty. These techniques will be discussed in this paper.
\end{abstract}

\begin{IEEEkeywords}
Unit Commitment, Stochastic Programming, Robust Optimization, Energy Grids.
\end{IEEEkeywords}

\section{Introduction}

Among the widespread area of mathematical optimization, the Unit Commitment (UC) offers a fundamental application on improvement issues wherever the assembly of a group of electrical generators is or has to be coordinated so as to realize a common output, sometimes either matching the energy demand at minimum price or maximizing revenue from electricity production \cite{b1}. This can be necessary as a result of it's troublesome to store current on a scale comparable to traditional consumption; hence, every variation within the consumption should be matched by a corresponding variation of the production \cite{b1}.

The short-time period operation of electrical generators is fraught with a huge variety of re-assets of uncertainty, consisting of network energy load, much less VER generation. In a state of affairs wherein a microgrid transacts strength on wholesale markets, the short-time period making plans can be exacerbated through the uncertainty related to power marketplace prices, as stated in \cite{b2}.

There is a developing frame of literature at the utility of stochastic programming  and likewise of Robust Optimization (RO) methods withinside the UC problem. As a relevant example, the authors in \cite{b2} examine the market price uncertainty based on a distributionally robust optimization (DRO) related to microgrids and leverage the Kullback-Leibler divergence to suggest a DRO solution for UC. Following the previous research, in \cite{b3} the authors apply the knowledge of the KL Divergence in function of a DRO-based algorithm in order to expressly assess the uncertainty associated the net load in absence of the probability distribution. Analogous to the preceding works, a study on greenhouse gas (GHG) is done in \cite{b3}, where visual representation of the recorded topics is presented, constraints and ex-ante evaluation of carbon tax payment are taken into account parameters; the result of \cite{b4} leads to a quantification and comparison of the one discussed in the research and the classical UC approaches.

\begin{table}[h]
	\caption{Nomenclature, variables and parameters as used in \cite{b1}}
	\begin{center}
		\begin{tabular}{c|r}
			\textbf{Symbol} & \textbf{Description} \\
			\hline
			$U$ & Set of feasable commitment decisions\\
			\hline
			$V$ &  Deterministic uncertainty set (range, region) \\
			\hline
			$u$ & Day-ahead decisions\\
			\hline
			$v$ & Uncertainty parameter\\
			\hline
			$c$ & Startup and shutdown costs \\
			\hline
			$q$ & Quadratic cost function (coefficient vector)\\
			\hline
			$\xi$ & Random Vector (Uncertainity)\\
			\hline
			$s$ & Realization\\
			\hline
			$p$ & Dispatch and Reserves\\
			\hline
			$f\left(\cdot\right)$ & Fuel costs\\
			\hline
			$F(u,v)$ & Real-time dispatch cost function\\
			\hline
			$A, B, H$ & Uncertain Matrices\\
			\hline
			$d$ & Uncertain demand vector\\
			\hline
			
		\end{tabular}
	\end{center}
\end{table}



In addition to this Introduction, the remainder of this paper goes as follows. Firstly, uncertainty as an input-source of the approach model for UC will be reviewed in Section II. Succeeding by this segment, an introduction to the various formulations and key algorithmic techniques that have been applied to UC will be given in section III. In Section IV, a study will be held on the advantages and disadvantages of the appointed solutions in Section III. Section V concludes the paper and considers possible future research directions.

\section{Uncertainity Modeling as an Input to UC}
Unit commitment problems such as addressing load, renewable power generation, and proper reaction to uncertainties, require an appropriate structure of uncertainty sets, which performs a vital position in figuring out the conservativeness of the model \cite{b5}. In this passage, we speak distinct tactics to assemble uncertainty units primarily based totally on ancient data, with the reason of decreasing the conservativeness while keeping the equal degree of robustness of the solution.

Multiple samples and increasing amounts of data on wind energy are more than well documented in historical research \cite{b6}--\cite{b7}. Electric power system operators are highly capable of satisfying demand, even when extreme events occur. These observations motivate the focus on risk avoidance, where risk can be reliably assessed as in SUC or RUC, rather than epistemic uncertainty about average performance as in other possible approach-models for UC.

The first method is to count on the net load for every time length at every node falls among a decrease and a higher bound \cite{b1}, which may be set with the aid of using positive percentiles of the random load and wind energy output primarily based totally on ancient data, as showed in \cite{b6}; or with the aid of using a set percent of the nominal load \cite{b6}. Note that a state of affairs specifies the net load for every hour and every unit inside the scheduling system. A parameter referred to as the price range of uncertainty is described to manipulate the deviation of all hundreds from their nominal values \cite{b1}.

Our second-stage technique is proposed in \cite{b7}, where in order to assemble uncertainty units, the usage of historic realizations of the random variables are connected through making use of convex units and a particular danger risk measures. To fulfill this connection, right here, the authors cognizance on secure inadequacy of power wind technology to satisfy the net load.

In the stochastic programming approach, the uncertainty is captured via discrete probabilistic eventualities, while withinside the robust optimization approach, the variety of its values is described via better-than-the-worst-case set \cite{b1}; in order to make an affordable comparison, it is essential to outline uncertainty units-sets from a regular source of eventualities derived from the stochastic programming model.

For RUC, building a right uncertainty set performs a critical position in figuring out the conservativeness of the model \cite{b1}--\cite{b4}. The uncertainty set is regularly described through a limited scope where parameters may not undercome a certain decrease; this applies likewise to the upper bound of the set; the reliability is primarily based on the suggested cost and volatility of the distribution \cite{b5}.

As a conclusion, the focus of robust optimization or stochastic programming with a risk measure, respectively, is on optimizing against the worst cases of an uncertainty set or distribution that can be estimated reliably \cite{b1}. Which approach is more appropriate depends on the decision maker's ability to use historical data to estimate descriptions of the uncertainty that will be valid for the study horizon, as well as on their level to elude possible risks taking into account worst cases \cite{b1}.

\section{UC Models under Uncertainity and Solutions Algorithms}

There is an ongoing and vigorous debate among the UC literature to evaluate or advise distinctive tactics for modeling the UC problem; the question remains on which technique affords the most efficacious outcomes.

Many distinctive formulations had been proposed to approach UC and lots of answer methodologies had been developed\cite{b1}--\cite{b10}. Over time those components and techniques have advanced and overcome their antecesors based on priority list and dynamic programming\cite{b1}, as new approaches of the past years the Kullback-Leibler Divergence\cite{b2} or Distributionally Robust Optimization\cite{b8} have been presented in order to satisfy modern-day maximum performance.

A quantity of measures were proposed to enhance grid operation and making plans with renewables, inclusive of balancing region consolidation, growing flexibility in the aid portfolio, demand-aspect management, and use of garage devices \cite{b1}.

To give an instance, Blanco and Morales propose in \cite{b8} an efficient Robust Solution to the two-stage stochastic unit commitment problem. Here was handled a reformulation of scenarios (inputs) for a two-stage UC problem under uncertainty that pursues an algorithm which can perform accordingly in both, ideal, and worst-case scenarios.

The proposed reformulation of \cite{b8} is primarily based totally on partitioning the pattern area of the uncertain elements through clustering the eventualities that approximate their chance distributions (A similar approach can also be seen in \cite{b6}--\cite{b7}). The degree of conservatism of the ensuing UC system is managed by the wide variety of partitions into which the stated pattern area is split.

To efficaciously resolve the proposed reformulation of the UC problem under uncertainty, the authors of \cite{b8} implemented a second parallelization and decomposition schemes that depend upon a so-mentioned column-and-constraint procedure. As a result, the numerical outcome shows that the proposed methodology is capable of dramatically lessen the specified computational time while enhancing the optimality of the answer found, so as the conclusion in \cite{b8}.

\subsection{Stochastic Programming}
Stochastic unit commitment (SUC) has been arising as a promising technology to handle with success power generation issues concerning uncertainties e.g., \cite{b1}. The concept of SUC is to make use of scenario-primarily based totally uncertainty illustration applying the UC formulation. Compared to truly the use of reserve constraints, stochastic methods have sure advantages, which includes cost, saving, and reliability improvement, as explain in \cite{b1}.

In stochastic programming, the concept of Conditional Values at Risk (CVaR) is widely employed as a tractable degree to visualize the hazard related to the imbalances of net load. In comparison to choice-limited models requiring extra binary variables, CVaR includes best linear constraints and non-stop variables, making it computationally attractive, as showed and proven in \cite{b1} and \cite{b10}.

Recalling to the research made in \cite{b8}, one can approach UC with help of the so-called: Two-Stage Models and Algorithms; this methodologies make proper use of two measures in terms of total expected cost. For this model, decisions will be divided in day-ahead against real-time choices; this can be mathematically expressed with the model taken from \cite{b1} as follows:

\begin{equation}
\underset{\textbf{u}\in U}{\min} \textbf{c}^T \textbf{u} + E_{\xi}\left[F(\textbf{u}, \xi)\right]
\end{equation}

\begin{equation}
	\begin{split}
		F(\textbf{u},s) &= min_{\textbf{p}_s, \textbf{f}_s} f\left(\textbf{p}_s\right) \\
		& A_s\textbf{u} + B_s\textbf{p}_s + H_s\textbf{f}_s \geq \textbf{d}_s
	\end{split}
\end{equation}

\subsection{Robust Optimization}
Besides stochastic and robust methods based on primarily uncertainty sets, Robust Optimization (RO) is a recently proposed method for UC where uncertain energy supplies are involved. In this method, the uncertain parameter is represented as a random variable-set whose distribution is composed by a lower bound which predicts the worst case scenario, from this point on all solutions happen to be feasible, increasing in accuracy\cite{b10}.

In assessment to modelling with stochastic programming techniques, Robust Unit Commitment (RUC) models try and include uncertainty disregarding the probability records from stochastic distributions, and as an alternative a range for the uncertainty is acquired\cite{b1}. As in the previously discussed methodology (SUC), were, searching for a minimum of the whole anticipated outcome is the main priority. Ruc centers its focus on minimizing the worst-case cost concerning all viable results of the uncertain factors. Assuredly a rather orthodox outcome can be expected with implementations of this sort; however, this can lead to preferable performance by deleting a wide range of non-feasible scenarios.\cite{b1}.

The current research on UC has pulled that RUC had been applied to cope with uncertainties specially from net load and electricity market price uncertainty\cite{b2}, microgrids and thermal generation resources \cite{b2}--\cite{b4}, or even security-constrained problems \cite{b9}.

A possible spectrum of solutions in this distribution can be achieved with the aid of using historical stochastic research as proposed in \cite{b1} and \cite{b10}, or a discrete set with a sure at the distances among possibility vectors as showed in \cite{b5}.

A mathematical model is likewise proposed by the authors of \cite{b1}, this go as follows:

\begin{equation}
\underset{\textbf{u}\in U}{\min} \left(\textbf{c}^T \textbf{u} + \underset{\text{v}\in V}{max} \left[F(\textbf{u}, \textbf{v})\right]\right)
\end{equation}

\begin{equation}
	\begin{split}
		F(\textbf{u}, \textbf{v}) &= \min_{\textbf{p},\textbf{f}} \textbf{q}^T\textbf{p}\\
		& A_{\textbf{v}}\textbf{u} + B_{\textbf{v}}\textbf{p} + H_{\textbf{v}}\textbf{f} \geq \textbf{d}_{\textbf{v}}
	\end{split}
\end{equation}

\section{Comparisons of Models and Algorithms}
Among diverse stochastic optimization strategies, Robust Optimization has been more sought through the grid operators \cite{b2}--\cite{b4}. Compared with the alternative strategies together with stochastic UC where problems concerning the variety of eventualities and computational complexity exist, Robust UC is much less computation overloaded and greater ideal to the device operators as their main operation precedence is to make sure the device reliability even withinside the worst-case scenario, that's assured through RUC \cite{b10}.

As surveyed on this paper, there has been a wealthy literature on the use of Stochastic and Robust Optimization techniques to enhance decision-making in the UC process. A qualitative comparison has been made so far between these two models: Robust optimization and stochastic dynamic programming. Quantitative comparisons are out of doors the scope of this paper, and may be observed in \cite{b4} or \cite{b9}. It is cited that even in the identical class of algorithms, distinct methods and corresponding algorithms might also additionally vary significantly \cite{b1}.

In summary, marketplace clearing with stochastic formulations is
a complicated issue. Revenue adequacy and related problems such
as pricing, settlement, marketplace power, and uplift expenses all need
To be addressed earlier than any version may be placed into use in practice \cite{b1}. Furthermore, the version have to be transparent, fair, and understandable for all of the marketplace contributors worried to attain consensus, as the research on papers \cite{b2} to \cite{b4}. It is anticipated that troubles together with the definition of eventualities, quantity of eventualities and definition of uncertainty units associated with SO might be pretty contentious, a number of the diverse stake-holders/regulators and continue to be a barrier to imposing SO withinside the marketplace clearing process \cite{b1}.

\section{Future Research}

\subsection{Uncertainty Modeling}
While the forecasting strategies for renewable era along with wind and sun have made giant development in current years, the error margin of such predictions posses nonetheless an meaningful power of value equivalent to the forecasting load in general \cite{b6}--\cite{b7}. Despite of the similarities of the forecast factors, the way to include prediction mistakes is important for all of the SO strategies mentioned above; the modeling and inclusion of uncertain mistakes strongly relies upon the accuracy of their parameters and at once ends in the selection of a selected SO framework so \cite{b5} (e.g., scenario-primarily based totally stochastic programming or uncertainty-set-primarily based robust optimization). The choice of a positive SO technique is likewise case based and noticeably prompted with the aid of using the hazard choice of the device operator. A sensitive stability desires to be performed among economics and reliability \cite{b10}. New studies on subjects in this region consist of the following.

\begin{itemize}
	\item Pioneering uncertainty modeling concepts. A full-size quantity of studies is being completed with inside-the-state-of-the-art research subjects that may be leveraged to cope with power generator problems \cite{b1}. For example, more contemporary subjects consisting of facts-pushed Robust optimization (RO) \cite{b9}, distributionally Robust optimization (DRO) \cite{b6}, and the Kullbach Langarian Optimization \cite{b3} promise in overcoming a few struggles for energetic grid applications. Distributionally Robust optimization assumes that the opportunity distribution of the uncertainty is not always properly known, and seeks to discover a hard and fast of cost-powerful answers that for all feasible opportunity distributions, are both continually possible or at the least possible with inside the worst case \cite{b6}. This is an evaluation of  stochastic models with risk measures, which normally expect the whole know-how of opportunity distributions of the underlying uncertainty.
	\item Better climate forecasting. As the goal of climate forecasting is pushed through how it's miles utilized in grid operations, renewable era forecasting and a way to comprise it in energetic operations are carefully coupled and need to be analyzed in an incorporated matter, as discussed in \cite{b4}, \cite{b6}, \cite{b7}.
	\item Improvement of current SO approaches. For instance, many questions along with the quantity of eventualities, possible scenarios reduction, and the assessment of the first-rate of eventualities in stochastic programming deserve nonetheless similarly, research. It is likewise really well worth investigating a way to combine the benefits of every character of stochastic strategies into grid operations. This proposal is widely proposed in \cite{b10}.
	\item Multiscale modeling. The more than one asset of uncertainties in device operations generally show off numerous temporal resolutions, starting from instant contingencies, to consistent call-for and renewable output fluctuations, to each day gas delivery variations. More designated modeling of various selections that correspond to the unfolding of the multiscale uncertainties will enhance version fidelity, and as a result assist deliver instructional techniques toward real-global implementation \cite{b5}.
\end{itemize}
\subsection{Computational Challenges and Approaches}
Advances in modeling additionally deliver new demanding situations to computational algorithms. In stochastic-programming, primarily based methodologies, actual operational practices and modeling necessities introduce discrete choices among the second scenario-layer, along with quick-begin generator rescheduling and threat constraints \cite{b9}. Efficient and powerful robustification strategies may be implemented to the SO to assist practice decomposition algorithms \cite{b10}.

For well-known SUC problems, acceleration strategies with extra powerful cuts may be incorporated in decomposition algorithms to accelerate the set of rules convergence rate \cite{b8}. When thinking about a totally big range of scenarios, green state of affairs, choice and discount strategies ought to be sought and incorporated with corresponding decomposition algorithms \cite{b8}.

In Robust-optimization-based solutions, superior algorithms for bilinear, or maybe trilinear programs, ought to be advanced to effectively resolve RUC hassle with extra well-known uncertainty sets\cite{b1}. Literature reformulations and algorithms primarily based totally on nested decomposition, like \cite{b8}, will be used to address adaptive RUC models with a couple of degrees and distributionally robust methods.  

\section{Conclusion}
Robust optimization and stochastic programming have been substantially mentioned and studied as options to optimize unit commitment schemes under uncertain conditions. A famous effect has arisen that the Robust approach, with its cognizance at the worst case, is higher capable of manipulate risk even as stochastic programming emphasizes predicted values \cite{b10}. However, the stochastic programming formula can without problems accommodate a danger measure. Moreover, the effects of each strategies rely strongly at the version for the unsure parameters, both the uncertainty set or the probabilistic eventualities hired withinside the optimization \cite{b10}.


\end{document}